%% file: agt-5-50.tex
\documentclass{gtart_h}
\input agtout

\lognumber{50}
\volumenumber{5}
\volumeyear{2005}
\papernumber{50}
\pagenumbers{1291}{1314}
\received{26 September 2003} 
\revised{16 September 2005}
\accepted{26 September 2005}
\published{6 October 2005}

\usepackage{graphicx,amssymb,psfrag}

\def\psfraga <#1,#2> #3#4{%
\psfrag {#3}{\smash{\rlap{\kern #1 \raise #2\hbox{#4}}}}}

\def\figref#1{\hyperlink{#1anchor}{Figure~\ref*{#1}}}
\def\anchor#1{\noindent\hypertarget{#1anchor}{\smash{$\phantom{99}$}}\newline}

\newcommand\mod{{\ \mathrm{mod}\ }}
\newcommand\Mod{{\mathrm{Mod}}}
\newcommand\tr{{\mathrm{tr}}}
\newcommand\Id{{\mathrm{Id}}}

\newcommand\inn{{\mathrm{in}}}

\newcommand\Hom{{\mathrm{Hom}}}

\newcommand\End{{\mathrm{End}}}

\newcommand\cen{{\mathrm{center}}}
\newcommand\sn{{\mathrm{sign}}}

\newcommand\st{{\ |\ }}

\newcommand\im{{i = 1, \ldots, m}}

\def\be{\begin{equation}}
\def\ee{\end{equation}}
\def\ba{\begin{array}}
\def\ea{\end{array}}
\def\bea{\begin{eqnarray}}
\def\eea{\end{eqnarray}}
\newcommand\cA{{\mathcal A}}

\newcommand\cC{{\mathcal C}}
\newcommand\cD{{\mathcal D}}
\newcommand\cE{{\mathcal E}}

\newcommand\cH{{\mathcal H}}

\newcommand\cN{{\mathcal N}}

\newcommand\cR{{\mathcal R}}

\newcommand\cT{{\mathcal T}}

\newcommand\cV{{\mathcal V}}

\newcommand\bU{{\mathbf{U}}}

\newcommand\BA{{\mathbb{A}}}

\newcommand\BC{{\mathbb{C}}}

\newcommand\BK{{\mathbb{K}}}

\newcommand\BN{{\mathbb{N}}}

\newcommand\BQ{{\mathbb{Q}}}
\newcommand\BR{{\mathbb{R}}}

\newcommand\BT{{\mathbb{T}}}

\newcommand\BZ{{\mathbb{Z}}}

\newcommand\fg{{\mathfrak{g}}}
\newcommand\fh{{\mathfrak{h}}}

\newcommand\fw{{\mathfrak{w}}}

\newcommand\fS{{\mathfrak{S}}}

\newcommand\ga{{\alpha}}

\newcommand\gd{{\delta}}

\newcommand\gl{{\lambda}}
\newcommand\gm{{\mu}}

\newcommand\gs{{\sigma}}

\newcommand\gY{{\Theta}}

\newcommand\gW{{\Omega}}

\newtheorem{thm}{{Theorem}}[section]
\newtheorem{lemma}[thm]{{Lemma}}
\newtheorem{prop}[thm]{{Proposition}}
\newtheorem{cor}[thm]{{Corollary}}

\theoremstyle{remark}
\newtheorem{rmk}{{Remark}}[section]
\newtheorem{eg}[rmk]{Example}
\newcommand{\thmref}[1]{Theorem~\ref{#1}}
\newcommand{\secref}[1]{Section~\ref{#1}}
\newcommand{\lemmaref}[1]{Lemma~\ref{#1}}
\newcommand{\corref}[1]{Corollary~\ref{#1}}
\newcommand{\propref}[1]{Proposition~\ref{#1}}

\newcommand{\rmkref}[1]{Remark~\ref{#1}}

\newcommand{\egref}[1]{Example~\ref{#1}}

\newcommand\il{{i = 1, \ldots, \ell}}
\newcommand\str{{\textrm{str}}}
\newcommand\fzx{{\BQ(\xi)}}
\newcommand\zx{{\BZ[\xi]}}
\newcommand\Ug{{\mathbf{U}(\fg)}}
\newcommand\UA{{\mathbf{U}\!_\cA}}

\newcommand\Uv{{\mathbf{U}_v}}
\newcommand\Uh{{\mathbf{U}_h}}
\newcommand\UK{{\mathbf{U}_\BK}}
\newcommand\Uzx{{\mathbf{U}_\zx}}
\newcommand\Ufzx{{\mathbf{U}_\fzx}}

\newcommand\TI{\BT}
\newcommand\VA{\cV\!_\cA}
\newcommand\VK{\cV_\BK}
\newcommand\AII{\BA}
\newcommand\AIIh{\BA[[h]]}
\newcommand\Vg{{\cV\!_\fg}}

\newcommand\oVl{{_1\! V_\gl}}
\newcommand\fVli{{_4\! V_{\gl_i}}}

\newcommand\tVmi{{_2\! V_{\gm_i}}}

\newcommand\tVl{{_2\! V_\gl}}
\newcommand\thVl{{_3\! V_\gl}}
\newcommand\fVl{{_4\! V_\gl}}
\newcommand\Zphi{Z}
\newcommand\hJ{\hat J}
\newcommand\Ttmp{T_2(\mu')}
\newcommand\Ttm{T_2(\mu)}

\begin{document}
\title{Almost integral TQFTs from simple Lie algebras}
\author{Qi Chen\\Thang Le}
\gtemail{\mailto{qichen@math.ohio-state.edu}, \mailto{letu@math.gatech.edu}}
\asciiemail{qichen@math.ohio-state.edu, letu@math.gatech.edu}

\asciiaddress{Department of Mathematics, The Ohio State 
University\\Columbus, OH 43210-1174, USA\\and\\School of 
Mathematics, Georgia Institute of Technology\\Atlanta, GA 30332-0160, USA}

\address{Department of Mathematics, The Ohio State 
University\\Columbus, OH 43210-1174, USA\\{\rm and}\\School of 
Mathematics, Georgia Institute of Technology\\Atlanta, GA 30332-0160, USA}

\begin{abstract}
Almost integral TQFTs were introduced by Gilmer \cite{g2}.  The aim of
this paper is to modify the TQFT of the category of extended
3-cobordisms given by \cite{tu1} to obtain almost integral TQFT.
\end{abstract}

\asciiabstract{Almost integral TQFTs were introduced by Gilmer [Duke
  Math. J. 125 (2004) 389--413].  The aim of this paper is to modify
  the TQFT of the category of extended 3-cobordisms given by Turaev (in
  his book: Quantum invariants of knots and 3-manifolds) to obtain an
  almost integral TQFT.}


\makeshorttitle

\section{Introduction}\label{intro}
Inspired by a 3-dimensional interpretation of the Jones polynomial for knots, Witten predicted 
that one can define topological invariants for 3-manifolds
using simple complex Lie algebras. The first concrete construction was obtained by Reshetikhin and 
Turaev \cite{rt} for $sl_2$. Soon after, similar invariants were constructed 
for all simple Lie algebras. They are called the Witten-Reshetikhin-Turaev 
invariants (WTR-invariants for short) or quantum invariants because Reshetikhin and Turaev's
construction is based on the theory of quantum groups.
It was also Witten's vision that the WTR-invariants can be extended to the 
Topological Quantum Field Theory (TQFT) in a rather natural manner.
Loosely speaking, a TQFT with base ring $K$ is a functor from the category $\cC$ to the category of $K$-modules where
the objects in $\cC$ are surfaces and morphisms in $\cC$ are 3-dimensional cobordisms.
We denote a TQFT by a pair $(\cT, \tau)$ where
$\cT$ and $\tau$ are maps between objects and morphisms respectively.
Most known TQFTs have anomalies. Anomalies can be resolved by introducing $p_1$-structure as in \cite{bhmv} or
by studying the category of extended 3-cobordisms as in \cite{tu1}. In this paper we follow the latter.
One of the features of extended 3-cobordisms is that they contain embedded colored ribbon graphs.
(See \secref{tqfts} for more details about TQFT and ribbon graph.)
For each simple complex Lie algebra $\fg$ and certain integer $r$, one can define a TQFT $(\cT^\fg_r, \tau^\fg_r)$.
In this case, the colors of the ribbon graphs come from the representations 
of the quantum group $\Uv(\fg)$ associated to $\fg$ with the parameter $v$ equal to a primitive
$r$-th root of unity.

Gilmer introduced the notion of almost integral TQFT in \cite{g2}. Let $D$ be a Dedekind domain contained in a ring $K$. 
A TQFT $(\cT, \tau)$ with base ring $K$ 
is called {\it almost ($D$-)integral} if there
exists some element $\cE \in K$ such that $\cE\tau(M)$ is in $D$ for any closed connected cobordism $M$.
(Note that in \cite{g2}, $\cE$ is required to be in $D$. One can show that most results in \cite{g2} hold in this more
general definition.) 
Gilmer showed that the TQFTs considered in \cite{bhmv} are almost integral. Using this,
Gilmer, Kania-Bartoszynska and Przytycki \cite{gkp} proved that for a prime $p>2$ if an 
integral homology 3-sphere is $p$-periodic, i.e.\ it admits a $\BZ /p \BZ$ action with fixed point 
set a circle, then its projective $sl_2$ (or $SO_3$) WRT-invariant satisfies some congruent relation.
This criterion is used to show that the Poincare sphere is 2, 3, and 5-periodic only in \cite{chenle1}.

In this paper, we show, in \thmref{thm_main}, that a modified version of the TQFT $(\cT_r^\fg, \tau_r^\fg)$
is almost integral. The modification is to restrict the colors for ribbon graphs embedded in extended 3-cobordisms.
The difficulty in proving \thmref{thm_main} is that one must show 
the integrality of the projective WRT-invariant for every closed 3-cobordism with a ribbon graph,
while it was sufficient in \cite{g2} to show the integrality 
for every closed 3-cobordism with a framed link. (Framed links are special ribbon graphs.)
\thmref{thm_main} is an important ingredient to prove that certain TQFTs with rings of algebraic integers as base rings can
be constructed from $(\cT_r^\fg, \tau_r^\fg)$, cf.\ \cite{chen}.
As a corollary of \thmref{thm_main}, the criterion in \cite{gkp} can be generalized to all simple Lie algebras (\corref{app2}).

We start by recalling necessary definitions in \secref{pre}. 
We state the main results in \secref{main}. 
We will discuss one application of \thmref{thm_main} in \secref{app}.
\secref{proofs} contains the proof of \thmref{thm_main}.    

The authors wish to thank the referee for numerous helpful comments.

\section{Preliminaries}\label{pre}
In this section, we recall some definitions and known results needed in this paper.

\subsection{Lie algebras and their quantum deformations}\label{lie}
We will define quantum group and then look at its representation theory.
The reader may skip this section if he or she is familiar with it.

Let $\fg$ be a simple complex Lie algebra with Cartan matrix ($a_{ij}$), $i, j = 1, \ldots, \ell$. Fix a Cartan subalgebra $\fh$ of $\fg$ and a set of simple roots $\Pi_\fh = \{\alpha_1, \ldots, \alpha_\ell\}$ in its dual space $\fh^*$. One can define a symmetric bilinear form ($\cdot | \cdot$) on $\fh^*$ in the following way. Multiply the $i$-th row of ($a_{ij}$) by $d_i \in \{1, 2, 3\}$ such that ($d_i a_{ij}$) is a symmetric matrix. Set ($\alpha_i | \alpha_j$) $= d_i a_{ij}$. This bilinear form is proportional to the dual of the Killing form restricted on $\fh$. Let $X$ and $Y$ be the weight lattice and the root lattice of $\fg$. The Weyl group $W$ acts on $X$ and $Y$ naturally. The order of the group $X/Y$ is det($a_{ij}$). Let $X_+$ be the set of dominant weights and $Y_+ = Y \cap X_+$. According to the general theory of Lie algebra, the finite-dimensional representations of $\fg$ are parameterized by the dominant weights. Let 
$d = \max_{1\leq i \leq l}\{d_i\}$. 
Then $d$ is the square of the ratio of the lengths of a long root and of a short root of $\fg$.
Let $\alpha_0$ be the highest short root, and let $\rho$ be half of the sum of positive roots. The dual Coxeter number is $h^\vee = 1 + \max_{\alpha > 0} (\alpha | \rho)/d$.

From now on we fix a simple Lie algebra $\fg$. Let $h$ be a formal parameter and 
$\Uh = \Uh(\fg)$ be the $\BC[[h]]$-algebra topologically generated by
the set of generators $\{X_i, Y_i, H_i\}_{1\leq i \leq l}$ and the relations
$$\ba{ll}
\, [H_i, H_j] = 0, &
[X_i, Y_j] = \delta_{ij} \frac{\sinh(hd_iH_i/2)}{\sinh(hd_i/2)}, \\
\, [H_i, X_j] =  a_{ij} X_j, & [H_i, Y_j] = -a_{ij} Y_j
\ea$$ 
and if $i \neq j$
$$
\sum_{k=0}^{1-a_{ij}} (-1)^k \left[ {1-a_{ij}}\atop{k}\right]_{v_i}
X_i^k X_j X_i^{1-a_{ij}-k} = 0
$$
$$
\sum_{k=0}^{1-a_{ij}} (-1)^k \left[ {1-a_{ij}} \atop {k}\right]_{v_i}
Y_i^k Y_j Y_i^{1-a_{ij}-k} = 0
$$
where $v_i = e^{h d_i/2}$ and 
$\left[ {a}\atop{b} \right]_x = 
[a]^!_x/[b]^!_x [a-b]^!_x$ for $a\geq b\geq 0$.
Here $[y]^!_x = \prod_{i=1}^{y} [i]_x$ 
and $[z]_x = (x^z - x^{-z})/(x - x^{-1})$. We use $[a]_i$ to denote $[a]_{v_i}$.

It's convenient to introduce some new elements in $\Uh$. Let
$$
v = e^{h/2},\quad \mathrm{so}\quad v_i = v^{d_i}
$$
and for $1 \leq i \leq \ell$
$$
K_i = v_i^{H_i} = e^{h d_i H_i/2}, 
$$
$$
E^{(n)}_i = E_i^n/[n]^!_i,\qquad F^{(n)}_i = F_i^n/[n]^!_i,
$$
where
$$
E_i = X_i K_i^{1/2}, \qquad F_i = K_i^{-1/2} Y_i.
$$
Suppose $\mu = \sum_i a_i \alpha_i \in Y$. Denote $\prod_i K_i^{a_i}$ by $K_\mu$.

One can put a Hopf algebra structure on $\Uh$ with coproduct $\Delta_h$, antipode $S_h$ and counit $\epsilon_h$ defined as follows.
$$
\Delta_h(H_i) = H_i \otimes 1 + 1 \otimes H_i,
$$
$$
\Delta_h(X_i) = X_i \otimes K_i^{1/2} + K_i^{-1/2} \otimes X_i,
$$
$$
\Delta_h(Y_i) = Y_i \otimes K_i^{1/2} + K_i^{-1/2} \otimes Y_i,
$$
$$
\epsilon_h(H_i) = \epsilon_h(X_i) = \epsilon_h(Y_i) = 0,
$$
$$
S_h(H_i) = -H_i, \quad S_h(X_i) = -K_i X_i, \quad S_h(Y_i) = -K_i^{-1} Y_i.
$$
Let $\BK$ be the field of fractions of $\BC[[h]]$ and $\UK = \Uh\otimes_{\BC[[h]]} \BK$. Let $\Uv$ be the $\BQ(v)$-subalgebra of $\UK$ generated by $\{E_i, F_i, K_i\}_{1\le i\le \ell}$. Let $\UA$ be the $\cA$-subalgebra of $\Uv$ generated by $\{E_i^{(n)}, F_i^{(n)}, K_i\, |\, 1\le i\le \ell, n\in\BN\}$ where  $\cA = \BZ [v, v^{-1}]$. The Hopf algebra structure on $\Uh$ induces Hopf algebra structures on $\UK$, $\Uv$ and $\UA$. See Chapter 3 of \cite{lu} for a proof of this fact. 

\begin{rmk}
Our $\Uh$ coincides with the one defined in \cite{kassel} except that our $v_i$ is denoted $q_i$ there and because of this our $\Uv$ is denoted $U_q$ there. We use the same coproduct as in \cite{kassel} which is opposite to the one used in \cite{lu}. Our $\Uv$, $\UA$ and $K_i$ are denoted $\bU$, $_\cA\bU$ and $\tilde K_i$ in \cite{lu}.
\end{rmk}

It is known that finite-dimensional simple modules of $\UK$ (resp.\ $\Uv$) are parameterized by the dominant weights of $\fg$, i.e.\ for any dominant weight $\gl$ there is a unique finite-dimensional simple module $\tilde V_\gl$ (resp.\ $V_\gl$) of highest weight $\gl$ and all finite-dimensional simple modules are of this form. Lusztig showed that $V_\gl$ has a (canonical) basis $B(\gl)$ such that the $\cA$-module, denoted $_{\cA}V_\gl$, generated by $B(\lambda)$ inherits a $\UA$-module structure from $V_\gl$. 

Let $\VK$ be the category of finite-dimensional $\UK$-modules of type I. The Hopf algebra structure on $\UK$ induces tensor product and duality in $\VK$. Any module $M$ from $\VK$ admits a weight space decomposition
\be\label{eq:decomp}
M = \bigoplus_{\gl\in X} M^\gl
\ee
where $M^\gl = \{x\in M\, |\, H_i(x) = (\ga_i | \gl)\, x\}$. Let $\VA$ be the category of finitely generated (as $\cA$-modules) $\UA$-modules which have weight space decomposition as in Equation (\ref{eq:decomp}). 
Then $_{\cA}V_\gl$ and its dual ${_{\cA}V_\gl}^* = \Hom_\cA(_{\cA}V_\gl, \cA)$ are objects 
of $\VA$ for $\gl \in X_+$. Clearly $\VA$ is also a tensor category with duality.

Let $\gY$ and $\bar \gY$ be the quasi-$R$-matrix and its inverse (Chapter 4 of \cite{lu}). They are infinite sums
\be\label{eq:quasi-R}
\gY = \sum_s a_s\otimes b_s, \quad \bar\gY = \sum_s \bar a_s\otimes \bar b_s
\ee
where $a_s, b_s, \bar a_s$ and $\bar b_s$ are in $\Uv$. So $\gY$ and $\bar \gY$ belong 
to some completion of $\Uv\otimes \Uv$. It turns out that they actually belong to some 
completion of $\UA\otimes \UA$, i.e.\ $a_s, b_s, \bar a_s$ and $\bar b_s$ are in $\UA$ 
(Corollary 24.1.6 of \cite{lu}). All the $a_s, b_s, \bar a_s$ and $\bar b_s$ act as 0 
on any $M\in \VK$ except for finitely many of them. Therefore it makes sense to consider 
the map $\gY$ (or $\bar \gY): M\otimes N \to M\otimes N$ for any $M, N$ in $\VK$.

For $M, N$ in $\VK$, define $\Psi : M\otimes N \to M\otimes N$ by
\be\label{eq:Psi}
\Psi(x\otimes y) = v^{(\nu | \mu)} x\otimes y,
\ee
if $x\in M^\nu$ and $y\in N^\mu$. Note that $(\nu | \mu)$ is not necessarily an integer\footnote{There exists an integer $\Delta$ for $\fg$ such that $(\nu | \mu)$ belongs to $\frac{1}{\Delta}\BZ$ if $\nu, \mu\in X$ and belongs to $\BZ$ if either $\nu$ or $\mu$ is in $Y$.}. Define the braiding operator $c_{M, N}$ as follows
\be\label{eq:Braiding}
c_{M, N} = P \Psi \bar\gY : M\otimes N \to N\otimes M
\ee
where $P(x\otimes y) = y\otimes x$. It commutes with the $\UK$ action and hence is an operator from $\VK$.

Let $\Omega = \sum_s S_h(a_s)b_s$ where $a_s$ and $b_s$ are the ones in Equation~(\ref{eq:quasi-R}). 
It is an element in some completion of $\UA$. For any $M\in \VK$ one can define $\Omega : M\to M$ 
because only finitely many terms in the summation act nontrivially. For such $M$ the twist 
operator\footnote{The twist operator is called the quantum Casimir operator in Chapter 6 of \cite{lu}.} is
\be\label{eq:Twist}
\theta_M : M \to M,\quad \mathrm{with} \quad \theta_M(x) = v^{(\nu + 2\rho | \nu)}\Omega(x),
\ee
for $x\in M^\nu$. The operator $\theta_M$ is invertible and commutes with $\UK$ actions. If $M = \tilde V_\lambda$, $\gl\in X_+$ then $\theta_M = v^{(\gl+2\rho | \gl)} \Id_M$.

The category $\VK$ is a ribbon category with braiding $c$ and twist $\theta$.
But $\VA$ is {\em not} because of the fractional powers in Equations (\ref{eq:Psi}) 
and (\ref{eq:Twist}). We will define a ribbon subcategory of $\VA$ in \secref{4cat}.

\subsection{TQFTs based on ribbon categories}\label{tqfts}
The objective of this section is to define TQFT of extended 3-cobordisms following Chapter IV of \cite{tu1} with some modification.
(See also Section 3.3 of \cite{bk}.) Manifolds are always orientable and smooth. Maps between manifolds are always smooth. Non-zero (tangent or normal) vectors are equivalent up to scalar multiple. We fix a ribbon category $\cR$ in the rest of this section.
As mentioned in the introduction, a key feature of extended 3-cobordisms is that they contain colored ribbon graphs, which we define next.

\subsubsection{Ribbon graphs}\label{rg}
A {\it band} is a homeomorphic image of $[0,1]$ with a non-zero normal field. The images of $0$ and $1$ are called the {\it bases} of the band. 
An {\it annulus} is a homeomorphic image of $S^1$ with a non-zero normal field. A band or an annulus is oriented if it is equipped with a non-zero tangent field. 
A {\it coupon} is a homeomorphic image of $[0,1]\times[0,1]$. The images of $[0,1]\times 0$ and $[0,1]\times 1$ are called the 
{\it bottom} and the {\it top} of the coupon respectively. Coupons are always oriented.

Let $M$ be a 3-manifold. A {\it ribbon graph} $\Omega$ in $M$ is a union of bands, annuli, and coupons embedded in $M$ such that:

(i)\qua The bands and the annuli are oriented.

(ii)\qua  $\Omega$ meets $\partial M$ at some bases of the bands (called free ends of $\Omega$). The normal vectors of the bands at these free ends are tangent to $\partial M$.

(iii)\qua  Other bases of bands (called fixed ends of $\Omega$) lie on coupons' top or bottom with normal vectors in the direction of the positive side of the coupon.

(iv)\qua  Bands, annuli, and coupons are disjoint otherwise. 

Let $\Omega$ be a ribbon graph. An {\it $\cR$-coloring} of $\Omega$ is an assignment to each band and annulus of $\Omega$ an arbitrary object of $\cR$, and to each coupon of $\Omega$ a morphism of $\cR$ in the following way. Let $C$ be a coupon in $\Omega$ and we are looking at its positive side with the top above the bottom. Suppose the fixed ends on the bottom (resp.\ top) of $C$ have colors, counting from left to right, $V_1, \ldots, V_m$ 
(resp.\ $W_1,\ldots, W_n$). Let $\epsilon_i$ be + (resp.\ $-$) if the $i$-th fixed end at the bottom of $C$ is oriented downward (resp.\ upward). Let $\gd_j$ be + (resp.\ $-$) if the $j$-th fixed end at the top of $C$ is oriented downward (resp.\ upward). We assign $C$ a morphism 
$$
f\in \Hom_\cR(\bigotimes_{i=1}^m V_i^{\epsilon_i}, \bigotimes_{j=1}^n W_j^{\gd_j})
$$
where we use the notation $V^+ = V$, $V^- = V^*$ for any object $V$ of $\cR$. A ribbon graph $\Omega$ together with an $\cR$-coloring $\gm$ is called an {\it $\cR$-colored ribbon graph}, denoted $\Omega(\gm)$. A ribbon graph is {\it partially $\cR$-colored} if some bands and/or annuli and/or coupons are $\cR$-colored.

\subsubsection{Extended surfaces}\label{es}
Another feature of extended 3-cobordisms is that their boundaries (or rather the boundaries of their underlying extended 3-manifolds, cf.\ \ref{e3m})
are extended surfaces.
An {\it $\cR$-extended surface} (or $e$-surface for short)
is a closed oriented surface $\Gamma$ together with a finite set of $\cR$-marks on it and a 
decomposable\footnote{For a closed surface $\Gamma$ with connected components $\Gamma_1,\ldots, \Gamma_m$,
a Lagrangian subspace of $H_1(\Gamma, \BQ)$ is {\it decomposable} if it is a direct sum of 
Lagrangian subspaces of $H_1(\Gamma_i, \BQ), \im$.} Lagrangian subspace of 
$H_1(\Gamma, \BQ)$. 
An {\it $\cR$-mark} on a closed surface $\Gamma$ is a point $p$ on $\Gamma$ associated
with a triple $(t, V, \nu)$, where $t$ ({\it direction} of the mark) is a non-zero 
tangent vector at $p$, $V$ ({\it color} of the mark) is an arbitrary object from $\cR$ 
and $\nu$ ({\it sign} of the mark) is $+$ or $-$. If $\Gamma$ is an $e$-surface we denote by $-\Gamma$ the $e$-surface
obtained from $\Gamma$ by reversing the orientation of $\Gamma$, keeping the Lagrangian subspace
and, for every $\cR$-mark on $\Gamma$, changing its sign while keeping its color and 
direction unchanged. The empty surface is considered as an $e$-surface which is not allowed 
to have any $\cR$-mark on it.

An {\it $e$-homeomorphism} between two $e$-surfaces is a homeomorphism between the underlying 
surfaces that respects the extended structure, i.e.\ 
orientation, $\cR$-marks and Lagrangian subspace.

\subsubsection{Extended 3-cobordisms}\label{e3m}
Before defining extended 3-cobordisms we need another notion. An {\it $\cR$-extended 3-manifold} is a triple $(M, \gW, w)$ that consists of an oriented 3-manifold $M$ with $\cR$-extended structure on $\partial M$, an integer $w$ ({\it weight} of $M$) and an $\cR$-colored ribbon graph $\gW$ (\secref{rg}) sitting in it. The boundary of $M$ is an $e$-surface such that: 

(i)\qua  The free ends of $\gW$ meet $\partial M$ at the $\cR$-marks only and each mark meets a free end of $\gW$.

(ii)\qua  The colors (resp.\ normal vectors) of the bands of $\gW$ and the colors (resp.\ directions) of the $\cR$-marks agree if incident. 

(iii)\qua  The sign of a mark is + (resp.\ $-$) if the incident free end is directed inward (resp.\ outward). 

(iv)\qua  The orientation on $\partial M$ is induced by that of $M$. 

(v)\qua  It is required that $w(\emptyset)=0$.

(There is no restriction on the Lagrangian subspace of $\partial M$.)

An {\it $\cR$-extended cobordism} (or $e$-cobordism for short) is a triple $(M, \Gamma, \Lambda)$ where $M$ is an
$\cR$-extended 3-manifold and $\partial M = (-\Gamma)\sqcup \Lambda$ is an $e$-surface. An {\it $e$-homeomorphism} 
between two $e$-cobordisms is a homeomorphism between the underlying 3-manifolds that respects the extended 
structures on 3-manifolds and their boundaries.

Suppose $(M, \Gamma, \Lambda)$ and $(M', \Gamma', \Lambda')$ are two $e$-cobordisms and there is an $e$-homeomorphism $f : \Lambda \to \Gamma'$. One can glue these two $e$-cobordisms along $f$ to get a new $e$-cobordism $(M\cup_f M', \Gamma, \Lambda')$. Here $M\cup_f M'$ is an $\cR$-extended 3-manifold with an $\cR$-colored ribbon graph (obtained by gluing gibbon graphs in $M$ and $M'$) sitting inside and weight computed as in IV.9.1 of \cite{tu1}.

\subsubsection{TQFTs based on $\cR$}\label{rtqft}
We finally can say what TQFTs really mean. Let $K$ be a commutative ring with unit and $\Mod(K)$ be the category of projective $K$-modules. 
A topological quantum field theory (TQFT) based on $\cR$ with ground 
ring $K$ is a pair $(\cT, \tau) = (\cT_\cR, \tau_\cR)$. Here $\cT$ is a modular 
functor\footnote{A functor $\cT$ is {\it modular} if $\cT(\Gamma\sqcup \Lambda)$ is naturally identified 
with $\cT(\Gamma)\otimes_K\cT(\Lambda)$ and $\cT(\emptyset) = K$, cf.\ III.1.2 of \cite{tu1}.}
(based on $\cR$ with ground ring $K$) from the category of $e$-surfaces and $e$-homeomorphisms to $\Mod(K)$
and $\tau$ assigns to every $e$-cobordism $(M, \Gamma, \Lambda)$ a $K$-homomorphism
$$
\tau(M) = \tau(M, \Gamma, \Lambda) : \cT(\Gamma) \to \cT(\Lambda)
$$
that satisfies the following four axioms.

(Naturality)\qua If $(M, \Gamma, \Lambda)$ and $(M', \Gamma', \Lambda')$ are two $e$-cobordisms and there is an $e$-homeomorphism $f : M \to M'$ then $\tau(M')\circ \cT(f|_\Gamma) = \cT(f|_\Lambda) \circ \tau(M)$.

(Multiplicativity)\qua $\tau(M\sqcup M') = \tau(M)\otimes_K \tau(M')$.

(Functoriality)\qua If $(M, \Gamma, \Lambda)$ and $(M', \Gamma', \Lambda')$ are two $e$-cobordisms and there is an $e$-homeomorphism $f : \Lambda \to \Gamma'$ then $\tau(M\cup_f M') = \tau(M') \circ \cT(f) \circ \tau(M)$. 

(Normalization)\qua Suppose that $\Gamma$ is an $e$-surface. Then $\tau(\Gamma\times [0, 1], \Gamma\times 0, \Gamma\times 1) = \Id_{\cT(\Gamma)}.$

Here $(\Gamma\times [0, 1], \Gamma\times 0, \Gamma\times 1)$ is the $e$-cobordism with weight 0 that is induced by the extended structure on $\Gamma$ in the obvious way. Let $(\cT, \tau)$ be a TQFT based on $\cR$ with ground ring $K$. For an $e$-cobordism $(M, \emptyset, \Gamma)$ one has $\tau(M): K \to \cT(\Gamma)$. By slight abuse of notation, we denote $\tau(M)(1)\in \cT(\Gamma)$ simply by $\tau(M)$.
The TQFT is called non-degenerate if for any $e$-surface $\Gamma$, the module $\cT(\Gamma)$ is generated over $K$ by the set
$$
\{\tau(M)\ |\ (M, \emptyset, \Gamma) \textrm{ is an $e$-cobordism}\}.
$$
A cobordism with a closed underlying 3-manifold is called closed. The image of a closed cobordism
$(M, \emptyset, \emptyset)$ under a TQFT is a linear map from $K$ to $K$. Therefore $\tau(M)=\tau(M)(1) \in K$. 
Let $D$ be a Dedekind domain contained in a ring $K$. Recall that A TQFT $(\cT, \tau)$ with base ring $K$ 
is called {\it almost ($D$-)integral} if there
exists some element $\cE \in K$ such that $\cE\tau(M)$ is in $D$ for any closed connected cobordism $(M, \emptyset, \emptyset)$.

\section{The main results}\label{main}
In this section we construct an almost integral TQFT based on a ribbon category coming from $\UA$-modules.

Recall that $\fg$ is a simple complex Lie algebra. In what follows we fix  an odd prime integer $r\ge d h^\vee$,
which does not divide $\det(a_{ij})$ in the notation of \secref{lie}. We also fix a Cartan subalgebra $\fh$
of $\fg$ and a set of simple roots $\Pi_\fh = \{\alpha_1, \ldots, \alpha_\ell\}$ in its dual space $\fh^*$.

The fundamental alcove of level $k$ is defined as 
$$
C_k = \{x \in \fh^* ~|~ (x|\alpha_i) \ge 0,\ (x|\alpha_0) \le k,\ \il\}
$$
where $k = r-1-(\rho | \alpha_0)$ and $\alpha_0$ is the short highest root associated to the simple roots we choose. Note that $k\ge 0$ because $r\ge dh^\vee$. ($C_k$ is equal to $\bar C_k$ in \cite{le2}.) The restriction that
$r\nmid\det(a_{ij})$ is to ensure 
the invertibility of the so-called $S$-matrix, cf.\ the paragraph before \egref{eg}.

\subsection{Four categories}\label{4cat}
As mentioned in the introduction, we need to restrict the colors for the ribbon graphs in extended 3-cobordisms. To describe these colors
we introduce the following 4 categories.

Let $\cV_1 = \cV_1(\fg, r)$ be a category of $\UA$-modules over $\cA$. Objects in $\cV_1$ are direct summands of $V_1\star\cdots\star V_n$ where $V_i$ is $_\cA V_{\lambda_i}$ or its dual ${_\cA V_{\lambda_i}}^* = \Hom_\cA(_\cA V_{\lambda_i}, \cA)$ with $\gl_i \in C_k\cap Y$ and $\star$ is either $\otimes_\cA$ or $\oplus$. Morphisms in $\cV_1$ are $\UA$-linear maps. So $\cV_1$ is a full subcategory of $\VA$. Obviously $\cV_1$ is closed under tensor product and duality. For any $\nu, \mu$ in the root lattice $Y$, $(\nu | \mu)$ is always an integer. The weights of an object from $\cV_1$ are in $Y$, and so the braiding and the twist can be defined as in Equations (\ref{eq:Braiding}) and (\ref{eq:Twist}) over the ring $\cA$. Hence $\cV_1$ is a ribbon category.

Let $\xi$ be a primitive $r$-th root of 1. Consider $\zx$ as an $\cA$-algebra by sending $v$ to $\xi$. Let $\Uzx = \UA\otimes_\cA\zx$. Let $\cV_2 = \cV_2(\fg, r)$ be a category of $\Uzx$-modules over $\zx$. The objects in $\cV_2$ are direct summands of $V\otimes_\cA\zx$ for $V$ in $\cV_1$. Morphisms in $\cV_2$ are $\Uzx$-morphisms. Obviously $\cV_2$ is closed under tensor product and duality. The braiding and twist for $\cV_1$ induce braiding and twist for $\cV_2$. Hence, $\cV_2$ is also a ribbon category.

Let $\Ufzx=\Uzx\otimes_\zx \fzx$ where $\fzx$ is the field of fractions of $\zx$. Let $\cV_3 = \cV_3(\fg, r)$ be a category of $\Ufzx$-modules over $\fzx$. The objects in $\cV_3$ are direct summands of $V\otimes_{\zx}\fzx$ for $V$ in $\cV_2$. Morphisms in $\cV_3$ are $\Ufzx$-morphisms. Similarly $\cV_3$ is also a ribbon category.

\begin{rmk}
The category $\cV_1$ is different from the category of finite dimensional $\fg$-modules. For example,
$_\cA V_{\lambda}^*$ is not isomorphic to $_\cA V_{-w_0(\lambda)}$ where $w_0$ is the longest element
in $W$, cf.\ Section 3 of \cite{chen}.
The categories $\cV_2$ and $\cV_3$ are isomorphic as braided tensor categories in general. (We do
not need this here. See \cite{chen} for a proof of this fact.)
They are quite complicated. For example, they are not semisimple, i.e.\ they contain some objects that are reducible but indecomposable.
To get a semisimple tensor category we need to quotient out some morphisms
in $\cV_3$. 
\end{rmk}

We follow Section 3 of \cite{kirillov} to define the last category. First note that objects in $\cV_3$ are tilting modules. An object $V$ in $\cV_3$ is called {\it negligible} if for every $f \in \End_{\cV_3}(V)$ we have $\tr_q(f) = 0$. Here $\tr_q(f) = \tr (K_{2\rho} f)$ is the quantum trace of $f$. Then the tensor product of an arbitrary object with a negligible object is negligible. The dual of a negligible object is also negligible. These can be seen easily using graphical calculus, cf.\ I.2.7 of \cite{tu1}.

Let $V$ and $W \in \cV_3$. A morphism $f : V \to W$ is negligible if $f = g g'$ for some $g' : V \to Z$ and $g : Z \to W$ such that $Z$ is negligible. Then for a negligible morphism $f:U\to U'$ its transpose $f^*$ is negligible. For any morphism $g:V\to V'$, $f\otimes g$ is negligible. If $U = V'$ (resp.\ $V=U'$) then $f g$ (resp.\ $g f$) is negligible. 
Hence it's possible to define the quotient category. (See also Section 3.3 of \cite{kt}.) Let $\cN = \cN(\fg, r)$ be the class\footnote{$\cN$ is an ideal of $\cV_3$ in the language of \cite{kt}.} of negligible morphisms in $\cV_3$. The quotient category $\cV_4 = \cV_4(\fg, r) = \cV_3/\cN$ has the same objects as $\cV_3$ and morphisms
$$
\Hom_{\cV_4}(V, W) = \Hom_{\cV_3}(V, W)/\textrm{negligible morphisms}.
$$
Any object $V$ in $\cV_1$ (resp.\ $\cV_2$, $\cV_3$) gives rise to $\bar V = V\otimes\zx$ (resp.\ $\hat V = V\otimes\fzx$, $[V] = V/\!\sim$) in $\cV_2$ (resp.\ $\cV_3$, $\cV_4$). Any morphism $f$ in $\cV_1$ (resp.\ $\cV_2$, $\cV_3$) gives rise to $\bar f = f\otimes 1$ (resp.\ $\hat f = f\otimes 1$, $[f] = f/\!\sim$) in $\cV_2$ (resp.\ $\cV_3$, $\cV_4$). Here $\sim$ is the equivalence up to negligible morphisms. We denote $_\cA{V}_\lambda$, $_\cA{\bar V}_\lambda$, $_\cA\hat{\bar V}_\lambda$ and $[_\cA\hat{\bar V}_\lambda]$ by $\oVl$, $\tVl$, $\thVl$ and $\fVl$ respectively.

Then $\cV_4$ is a modular category that is dominated by simple objects $\fVl$, $\lambda \in C_k\cap Y$
(cf.\ Remark 3.10 of \cite{kirillov}\footnote{The category $\cV_3$, which is enough to produce a modular category,
is a category of tilting modules contained in $Rep\ U_{\xi}$ in \cite{kirillov}. Our $d$ is denoted $m$ there and
$\varkappa = r/2d$.}). See \cite{tu1} for the definition of modular categories. The invertibility of the so-called $S$
matrix in our case is proved in Theorem 3.3 of \cite{le2}, where $v$ is set to be a primitive $2r$th 
root of 1. But this causes no problem. The $S$ matrix (before plugging in $\xi$) is a matrix over $\BZ[v^2, v^{-2}]$ so in either case the entries of the $S$-matrix are substitutions of $v^2$ by a primitive $r$th root of 1. (Remember
that $r$ is odd.)

\begin{eg}\label{eg}
Let $\fg = sl_2$. The objects $\fVl$,
$\lambda = 0, \alpha_1, \ldots, \frac{r-3}{2} \alpha_1$, are irreducible in $\cV_4$ and 
every object in $\cV_4$ is a direct sum of them. On the other hand, ${_2\! V_{\alpha_1}}\otimes{_2\! V_{\frac{r-1}{2}\alpha_1}}$ is not a direct sum of simple objects in $\cV_2$.
\end{eg}

\subsection{The TQFT}\label{tqft}
In this section we define a TQFT, based on the ribbon category $\cV_2$ with ground ring $\fzx$ (cf.\ \secref{rtqft}),
from the category of $\cV_2$-extended cobordisms to $\cV_4$. 
We modify the construction given in IV of \cite{tu1}. 
In the rest of this paper $e$-surfaces and $e$-cobordisms will mean 
$\cV_2$-extended surfaces and $\cV_2$-extended cobordisms.

\subsubsection{Parameterization of $e$-surfaces}\label{pe}
\begin{figure}\anchor{hb}
\centering
\psfrag{1}{
}
\psfrag{2}{
}
\psfrag{3}{
}
\psfrag{4}{
}
\includegraphics[width=4in]{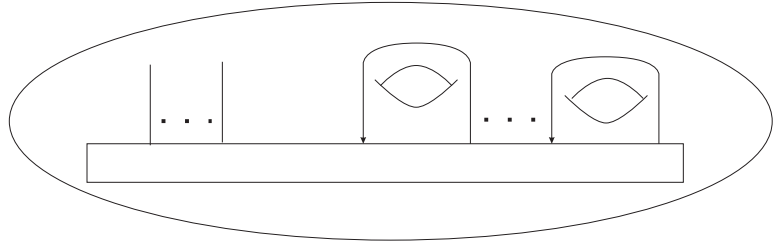}
\caption{A standard handlebody 
}\label{hb}
\end{figure}

We also need to parametrize every $e$-surface by a {\it standard surface}, 
which by definition is either an empty surface or an $e$-surface that is the boundary of a standard
handlebody (with induced extended structure). Here by 
a {\it standard handlebody} we mean a genus $g$ handlebody standardly embedded
in $\BR^3$ with a partially $\cV_2$-colored oriented ribbon 
graph $R$ sitting inside as shown in \figref{hb}.
See \secref{rg} for the definition of ribbon graphs.
The ribbon graph $R$ consists of a coupon (the narrow rectangle near the bottom), several vertical bands 
(with fixed ends on the coupon and free ends on the boundary of the standard handlebody) and $g$ half-circled 
bands (oriented to the left with bases on the coupon). The vertical bands are oriented and colored by
objects from $\cV_2$. The coupon and the $g$ half-circled bands are not colored.

For any connected $e$-surface $\Gamma$ let 
$$
\{p: \Sigma \to \Gamma \st \Sigma \textrm{ is a standard surface and } 
p \textrm{ is an $e$-homeomorphism}\}
$$
be the set of all parameterizations of $\Gamma$ up to $e$-isotopy (in the obvious sense). 
For any $e$-surface $\Gamma$, this set is not empty (IV.6.4.2 of \cite{tu1}).
Note that $\Gamma$ may be parameterized by more than one standard surfaces because the 
$\cV_2$-marks on $\Gamma$ are not ordered.

\subsubsection{The modular functor}\label{mf}
Let's define a modular functor $\cT$.
First, let $\cT(\emptyset) = \fzx$.

Suppose $\Gamma$ is a genus $g$ 
$e$-surface with $m$ $\cV_2$-marks $(t_i, V_i, \epsilon_i)$, $1\le i\le m$ (cf.\
\secref{es}). For each parameterization $p$ of $\Gamma$, put
\begin{eqnarray}
	\cT_p(\Gamma) = \bigoplus_\gl \Hom_{\cV_4} \left(\fzx, \bigotimes_{i=1}^m [\hat V_i]^{\epsilon_i}\otimes
	\bigotimes_{i=1}^g(\fVli\otimes \fVli^*)\right)
\end{eqnarray}
where $\gl = (\gl_1,\ldots,\gl_g)$ ranges over $(C_k\cap Y)^g$.
Recall from \secref{4cat} that any object $V$ in $\cV_2$ induces an object $[\hat{V}]$ in $\cV_4$.
Again $[\hat{V}_i]^+ = [\hat{V}_i]$ and $[\hat{V}_i]^- = [\hat{V}_i]^*$.

Furthermore, for two parameterizations $p, p'$ of $\Gamma$ there is an isomorphism 
$\varphi(p, p') : \cT_p(\Gamma) \to \cT_{p'}(\Gamma)$ of $\fzx$-vector 
spaces, defined as in IV.6.3 and IV.6.4.2 of \cite{tu1}. 
Identify the vector spaces 
$$
\{\cT_p(\Gamma) \st p \textrm{ is a parameterization of } \Gamma\}
$$ 
along the isomorphisms 
$\{\varphi(p, p')\}_{(p,p')}$. The resulting vector spaces $\cT(\Gamma)$ depends only on 
$\Gamma$. For any parameterization $p$ of $\Gamma$, $\cT(\Gamma)$ is canonically 
isomorphic to $\cT_p(\Gamma)$. Denote this isomorphism by $p_\sharp$.




Let $\Gamma$ and $\Gamma'$ be connected $e$-surfaces. For an $e$-homeomorphism 
$f: \Gamma\to \Gamma'$ we define $\cT(f): \cT(\Gamma)\to \cT(\Gamma')$ as follows. 
Pick any parameterization $p: \Sigma \to \Gamma$. Then $f p$ is a parameterization 
of $\Gamma'$. Set $\cT(f) = (f p)_\sharp (p_\sharp)^{-1}$ which does not depend on 
the chosen parameterization (cf.\ IV.6.3.1 of \cite{tu1}).

For non-connected $e$-surfaces we can do the above componentwisely and then form the tensor product. 
The above process defines a modular functor $\cT$ from the category of $e$-surfaces and $e$-homeomorphisms to
$\Mod(\fzx)$ (cf.\ IV.6.3.3 of \cite{tu1}).


\begin{rmk}
Unlike \cite{tu1}, which uses only a modular category in the construction, 
we start with the ribbon category $\cV_2$ in the construction of the $e$-surfaces 
and $e$-cobordisms and use the modular category $\cV_4$ later when we define the modular 
functor $\cT$. This modification ensures that the TQFT constructed in the next subsection is 
almost integral. See also \rmkref{free}.
\end{rmk}

\begin{rmk}
We quote several results from \cite{tu1} in this and the next subsections. They are valid in our construction because the places where we use these results are the places where we are using the modular category $\cV_4$.
\end{rmk}

\subsubsection{An almost integral TQFT}\label{anitqft}
Now we can describe the almost integral TQFT $(\cT, \tau)$ as claimed at the beginning of \secref{main}. Let $(\cT^{\cV_4}, \tau^{\cV_4})$ be the TQFT
derived from the modular category $\cV_4$
in the way given in Section IV.9 of \cite{tu1}.
By definition, $\cT^{\cV_4}(\check\Gamma) = \cT(\Gamma)$ where $\check\Gamma$ denotes the $\cV_4$-extended surface induced from $\Sigma$ by applying $[\hat{\ }]$ (cf.\ \secref{4cat})
to the colors of the $\cV_2$-marks on $\Gamma$.
For a $\cV_2$-extended cobordism $M$,
we define $\tau(M)$ to be $\tau^{\cV_4}(\check M)$,
where $\check M$ denotes the $\cV_4$-extended cobordism induced from $M$ by applying $[\hat{\ }]$
to the colors of the $\cV_2$-extended ribbon graph in $M$ and the colors of the $\cV_2$-marks on $\partial M$.
   

\begin{lemma}\label{non-degenerate}
The pair $(\cT, \tau)$ is a non-degenerate TQFT, based on $\cV_2$ with ground ring $\fzx$,
from the category of $\cV_2$-extended cobordisms to $\cV_4$.
\end{lemma}
\begin{proof}
The claim that $(\cT, \tau)$ is a TQFT is almost immediate from the definition because we actually define the pair through a TQFT $(\cT^{\cV_4}, \tau^{\cV_4})$. We indicate why it's non-degenerate, i.e.\ $\cT(\Gamma)$ is generated by $\tau(M)$ with $\partial M = \Gamma$. (Recall that we denote $\tau(M)(1)$ by $\tau(M)$.) It suffices to consider the case when $\Gamma$ is connected and $\Gamma = \Sigma$ for some standard surface $\Sigma$.

Recall that $\Sigma$ bounds a standard handlebody $\cH$ with a partially $\cV_2$-colored
ribbon graph $R$ sitting inside. Let $\cH'$ be the standard handlebody obtained from $\cH$ by forgetting the $\cV_2$-coloring on $R$. We denote the uncolored ribbon graph in $\cH'$ by $R'$. Suppose the vertical bands in $R$ are colored by $V_1,\ldots, V_m$. By an {\it admissible $\cV_4$-coloring} (or $a$-coloring for short) of $\cH'$ we mean a $\cV_4$-coloring of $\cH'$ such that the $j$-th vertical band of $R'$ is colored by $[\hat{V}_j]$ and the $i$-th half-circled band is colored by $\fVli$ with $\gl_i\in C_k\cap Y$. 
According to IV.2.1.3 of \cite{tu1} the set $\{ \tau^{\cV_4}(\cH'(\mu))\}_\mu$ generates $\cT^{\cV_4}(\check \Sigma) = \cT(\Sigma)$ where $\mu$ runs through all $a$-colorings of $\cH'$. We have to generate this set over $\fzx$ using $\cV_2$-extended cobordisms. 

Any morphism $f\in \cV_4$ can be lifted (non-uniquely) to a morphism 
$f'\in \cV_3$, for which there exists $A\in \zx, A\neq 0$ such that $A f' = \hat{f''}$ for some 
$f''\in \cV_2$. Any $a$-coloring $\mu$ of $\cH'$ is uniquely determined by its color on the coupon,
denoted $f_\mu$. We can find a $f_\mu''\in \cV_2$ as above. Let $\mu''$ be the 
$\cV_2$ coloring of $\cH'$ such that the coupon is colored by $f''_\mu$, the 
vertical bands are colored by $V_j$ and the half-circled bands are colored 
according to $f''_\mu$. Then $\cH'(\mu'')$ is a $\cV_2$-colored cobordism. 
It's clear that $\tau(\cH'(\mu'')) = A \tau^{\cV_4}(\cH'(\mu))$ for some 
$A \in \zx, A\neq 0$. (The nonuniqueness of the lifts causes no trouble 
because the difference between two lifts of one morphism is negligible.)
\end{proof}

\begin{rmk}\label{free}
We did not start with $\cV_4$ because that would give too much freedom in coloring
ribbon graphs and that in turn would destroy any integrality of the TQFT. 
Here is why. Recall that $(\cT^{\cV_4}, \tau^{\cV_4})$ is the TQFT based on $\cV_4$ with ground ring 
$\fzx$ constructed as in \cite{tu1}. We claim that it can not be almost integral.
Suppose $D$ is a Dedekind domain in $\fzx$. Let $M$ be a closed connected $\cV_4$-extended 
cobordism which contains a coupon colored by $f$. Let $M_a$ be the $\cV_4$-extended
cobordism obtained from $M$ by changing $f$ to $a f$ for some $a\in\fzx$. It's impossible 
to find a universal constant $\cE\in \fzx$ such that $\cE \tau(M_a)\in D$ for all $a\in\fzx$. 
\end{rmk}

Recall that $\ell$ is the rank of $\fg$ and $w_0$ is the longest element in its Weyl group. Let $\zeta$ be a root of unity such that
$$
\textrm{order}(\zeta) = \left\{\ba{ll}
r & \textrm{if\ $\ell$\ is\ even,\ $r$\ is\ arbitrary\ and\ }\ \sn(w_0)=1,\\
r & \textrm{if\ $\ell$\ is\ odd\ and\ } \sn(w_0)r\equiv 1 \mod 4,\\
4r & \textrm{otherwise.}
\ea\right.
$$
The following theorem is our main result which will be proved in \secref{proofs}.

\begin{thm}\label{thm_main}
Let $r\ge d h^\vee$ be an odd prime that is not a factor of $\det(a_{ij})$.
Then the TQFT $(\cT, \tau)$ is almost $\BZ[\zeta]$-integral.
\end{thm}

\subsection{Calculating $\tau$}\label{sec:tau}
Let $(M, \gW, w)$ be a closed connected $\cV_2$-extended 3-manifold. We recall how to calculate $\tau(M)$.
Let $L$ be a framed link in $S^3$ such that $M$ is obtained from $S^3$ by surgery along $L$. 
Then there exists a $\cV_2$-colored ribbon graph 
$\gW'$ in $S^3$ disjoint from $L$ such that $(M, \gW)$ is obtained from $(S^3, \gW')$ by 
surgery along $L$. By slight abuse of notation we will write $\gW$ for $\gW'$. 

Suppose that $L$ has $m$ components. Let $\gm = (\gm_1, \ldots, \gm_m) \in (C_k\cap Y)^m$.
Denote by $L(\gm)$ the $\cV_2$-colored framed link with the $i$-th component colored by $\tVmi$. Let $U^m$ be the trivial link of $m$ components. Let $J = J^\fg$ be the quantum invariant of $\cV_2$-colored ribbon graphs in $S^3$ defined in I.2.5 of \cite{tu1}. Define
\be\label{F}
F_{(L, \gW)} = \sum_{\gm\in (C_k\cap Y)^m} J_{U^m(\mu)} J_{\gW\sqcup L(\mu)}.
\ee
We put $F_\pm = F_{(U_\pm, \emptyset)}$ for the trivial knot $U_\pm$ with $\pm 1$ framing.
Let $\kappa$ and $\eta$ be one of the 
square roots of $F_-/F_+$ and $F_- F_+$ respectively such that $\kappa \eta = F_-$. (Note that $\eta$ 
and $F_-$ are denoted $\cD$ and $\Delta$ in \cite{tu1}. Hence our $\kappa$ is $\Delta/\cD$.)
According to IV.9.2 and II.2.2 of \cite{tu1}
\bea\label{tau}
\tau(M) & = & \tau (M, \gW, w)(1) = F_-^\gs 
F_{(L, \gW)}\eta^{-(\gs + m +1)} \kappa^w \\
& = &  F_{(L, \Omega)}
\eta^{-(m+1)} \kappa^{w+\sigma} \nonumber\\
& = & \frac{F_{(L, \Omega)}}{F_-^{\sigma_- + \beta_1}
F_+^{\sigma_+}}
\eta^{-1} \kappa^{\beta_1 + w} \nonumber
\eea
where $\beta_1$ is the first Betti number of $M$ and $\gs, \gs_+$ and $\gs_-$ are the signature, the number of positive eigenvalues and the number of negative eigenvalues of the linking matrix of $L$ respectively. Recall that $\sigma = \sigma_+ - \sigma_-$ and $m = \sigma_+ + \sigma_- + \beta_1$.

\section{One application}\label{app}
Again let $r\ge d h^\vee$ be an odd prime that is not a factor of $\det(a_{ij})$ and let $\xi$ be a primitive $r$-th root of 1. It is well known that for a closed connected oriented 3-manifold $M$, the projective WRT-invariant $\tau^{P\fg}_M(\xi)$ is equal to $\eta \tau(M)$ multiplied by some algebraic integer for a universal complex number $\eta$, cf.\ Equation~(\ref{tau}).
Here $M$ is considered as a closed $\cV_2$-extended cobordism with 0 
weight and empty ribbon graph sitting inside. 
It is also known that $\tau^{P\fg}_M(\xi)$ is in $\zx$.

One classical question in low dimensional topology is to determine whether a 3-manifold
is $p$-periodic. A manifold is $p$-periodic if it admits $\BZ/p\BZ$ action with fixed point 
set a circle. It is shown in \cite{gkp} that
if a homology sphere $M$ admits such an action then $\tau^{Psl_2}_M$ will satisfy some congruent
relation. One can generalize it to all quantum invariants using \thmref{thm_main}.

\begin{cor}{\rm\cite{chenle1}}\label{app2}\qua
Let $M$ be an $r$-periodic integral homology 3-sphere, and let $r$, $\xi$ be as above. Then
$$
\tau^{P\fg}_M(\xi) \equiv \xi^s \overline{\tau^{P\fg}_M(\xi)} ~\mod~ 
r ~\inn~ \zx
$$
for some integer $s$. Here `bar' is the complex conjugation.
\end{cor}

\begin{rmk}
Let $P$ be the Poincare sphere ($-1$ surgery on the left-hand trefoil), and let $B$ be the Brieskorn sphere ($-1$ surgery on the right-hand trefoil). Corollary \ref{app2} is used in \cite{chenle1} to show that $P$ has periodicity 2, 3, 5 only and $B$ has periodicity 2, 3, 7 only.
\end{rmk}

\section{Proofs}\label{proofs}
\subsection{Proof of \thmref{thm_main}}\label{proof_main}
In this section we  use \propref{prop_poly}, whose proof will be postponed till \secref{poly}, to prove \thmref{thm_main}.  Recall that $r$ is an odd prime $\ge d h^\vee$ which 
does not divide $\det(a_{ij})$ and $\xi$ is a primitive $r$th root of 1.
Let $\gW\sqcup L$ be a ribbon graph in $S^3$ where $\gW$ is a $\cV_2$-colored ribbon graph
(possibly with some $\cV_2$-colored annulus components) 
and $L$ is an uncolored framed link of $m$ components.

\begin{prop}\label{prop_poly}
With the above notation, $F_{(L, \gW)}$ is in $\zx$ (cf.\ Equation~(\ref{F})), and
it is divisible by $(\xi - 1)^{m(r\ell - \dim \fg)/2}$ in $\zx$.
\end{prop}

This proposition will be proved in \secref{poly}.

\begin{proof}[Proof of \thmref{thm_main}]
Let $(M, \gW, w)$ be a closed connected $\cV_2$-extended 3-manifold. Suppose that $M$ is the result 
of surgery along a framed link $L \subset S^3$ of $m$ components. By Equation~(\ref{tau}), 
$$
\tau(M) = \frac{F_{(L, \Omega)}}{F_-^{\sigma_- + \beta_1} F_+^{\sigma_+}} \eta^{-1} \kappa^{\beta_1 + w}.
$$
By Proposition 4.4 in \cite{le2},
$
F_\pm / (\xi - 1)^{(r\ell-\dim \fg)/2}
$
is invertible in $\zx$.
Hence, using \propref{prop_poly}, $\frac{F_{(L, \Omega)}}{F_-^{\sigma_- + \beta_1} F_+^{\sigma_+}}$ is in $\zx$.
This implies that $\eta \tau(M) \in \BZ[\xi, \kappa]$. The theorem now follows from the next lemma.
\end{proof}

\begin{lemma}
With the notations as above, $\kappa\in \BZ[\zeta]$.
\end{lemma}

\begin{proof}[Sketch of proof]
First, one can follow the proof of Theorem 6.2 of \cite{kirillov} to show that 
$\eta^2 = F_+ F_- = \sn(w_0)r^\ell/a^2$ for some $a\in\zx$. Second, it's known that $\sqrt{\sn(w_0) r^\ell} \in \BZ[\zeta]$.
Hence $\eta$ belongs to $\BQ(\zeta)$. Since $\eta\kappa = F_-$ is in $\zx$, $\eta$ is in $\BQ(\zeta)$ too.
It remains to notice that $\kappa$ is always a root of unity (cf.\ again \cite{kirillov}).
\end{proof}

\subsection{Proof of \propref{prop_poly}}\label{poly}
The proof of \propref{prop_poly} relies on the Kontsevich integral for framed tangles (cf.\ \cite{lm, kt}).
The Kontsevich integral can be considered as a functor from the category of framed tangles to the category of chord diagrams.
We define these two categories first.

\subsubsection{Categories of framed tangles and chord diagrams}\label{tangle}
Up to the end of \secref{weight}, 
we follow \cite{kt} closely so that the reader can easily find some missing details there. 
For a different approach see \cite{lm}.

Let $\BT$ be the category of framed tangles. The objects are words in +'s and $-$'s and the morphisms are framed tangles (up to the equivalence classes) in $\BR^2\times [0, 1]$. The category $\BT$ is a strict ribbon category (cf.\ Section 2.1 of \cite{kt}, where $\BT$ is denoted $\cT$).

For any pair of non-negative integers $(a, b)$, an {\it $(a, b)$-curve} $C$ 
is a compact oriented smooth curve\footnote{Smooth curves are smooth 1-manifolds. They do not have to be embedded in 3-manifolds nor have associated normal vectors.} such that $\partial C$
is divided into two totally ordered sets, source $s(C)$ and target $t(C)$ 
of cardinalities $a$ and $b$ respectively.
An {\it $(a,b)$-chord diagram} is an $(a, b)$-curve $C$ together with finitely many pairs
of unordered points (called chord-ends)
on $C$. Such a pair is usually indicated 
by a dashed line (called a chord on $C$)
connecting the pair. Let $\BA$ be the category
of chord diagrams. An object (resp.\ a morphism) in $\BA$ is a word in +'s
and $-$'s (resp.\ a $\BC$-linear combination of chord 
diagrams (up to the equivalence classes) with the 
same target and source). Then $\BA$ is a strict infinitesimal symmetric category.
In Section 2.3 of \cite{kt}, $\BA$ is denoted $\cA(\BC)$. Note that
chord diagrams in \cite{kt} are equipped with residues, an extra structure
being omitted here for it is not essential for our purpose.

\subsubsection{The Kontsevich integral for framed tangles}\label{kon}
For any strict infinitesimal symmetric category one can obtain a braided tensor category using formal integral, cf.\
the proof of Theorem 4.7 of \cite{kt}. We recall it when the strict infinitesimal symmetric category is $\BA$.

Let $h$ be a formal parameter. The {\it formal integral} $\AIIh$ of $\AII$ is a category whose objects are the 
same as those of $\AII$ and whose morphisms from $s$ to $s'$ are formal series $\sum_{i=0}^\infty f_i h^i$ with 
$f_i \in \Hom_{\AII}(s, s')$. There exists a braided tensor category structure on
$\AIIh$. This category is {\it not} a strict braided tensor category.

Let $\AII'=\AIIh^{\str}$ be the strictified category of $\AIIh$. See XI.5 of \cite{kassel} for the definition
of strictification. Then $\AII'$ is a strict braided tensor category
whose objects are finite sequences of objects of $\AII$, including the empty sequence.
For such a sequence $s = (s_1, \ldots, s_k)$, let $P(s) = \otimes_i s_i$ if $k>0$, where tensor product is taken in the order of index $i$, and $P(s) = \BC$ if $k=0$. 
For any two objects $s$ and $s'$ in $\AII'$ let
$$
\Hom_{\AII'}(s, s') = \Hom_{\AIIh}(P(s), P(s')).
$$

\begin{lemma}{\rm\cite{kt} Section 6}\label{lemma_kon}\qua
There is a functor 
$\Zphi : \TI \to \AII'$ respecting braided tensor category structure.
\end{lemma}

The functor $Z$ is called the {\it Kontsevich integral for framed tangles}. For any framed tangle $T$, 
\be\label{eq:kon}
Z(T) = \sum_{i=0}^\infty Z_i(T)h^i
\ee
where $Z_i(T)$ is a linear combination of chord diagrams with $i$ chords,
whose underlying 1-manifolds are $T$ forgetting framing.
We call $Z_i(T)$ the degree $i$ part of $Z(T)$.

\subsubsection{A weight system}\label{weight}
To calculate the quantum invariant $J^\fg$ from the Kontsevich integral we need to define a weight system.
We are only interested in the string tangle case.
A {\it string tangle} is a tangle whose $i$-th free
upper end is connected to the $i$-th free lower end for all $i$. A string
tangle is allowed to have closed components. A {\it string chord diagram} is defined similarly.

Let $V$ be a finite dimensional $\fg$-module.
Let 
$
C : \epsilon_1\cdots\epsilon_m
\to \epsilon_1\cdots\epsilon_m
$
be a string chord diagram. It is a morphism in $\AII$.
A state $s_0$ of $C$ is an assignment to each unordered pair of points (chord) 
an element $x_i$, $i \in \{1, \ldots, \dim \fg\}$ where $\{x_i\}$ is an orthonormal basis of $\fg$ with respect to the Killing form. For each component of the underlying 1-manifold of $C$ we get an element in the universal enveloping algebra $\Ug$ obtained by multiplying
$x_i$'s on this component in the order of the orientation of this component. Thus we get
$n$ elements $y_j$ in $\Ug$ where $n$ is the number of components of $C$.
Tensoring up these $n$ elements one gets an element in $\Ug^{\otimes n}$.
If the $i$-th component of $C$ is an annulus then we take the trace 
$\tr_V(y_i)$. Having done so for all annulus components we end up with
$w(s_0) \in \Ug^{\otimes m}$. The weight of $C$ is
\be\label{eq:we}
\fw_\fg(C) = \sum_{s_0} w(s_0) : \bigotimes_{i=1}^m V^{\epsilon_i} \to \bigotimes_{i=1}^m V^{\epsilon_i}
\ee
summed over all possible states of $C$, where $V^+=V$ and $V^-=V^*$.
Note that although $y_i$ is not well-defined for an annulus component, the
weight is well-defined because we take the trace. We extend the weight to
linear combination of chord diagrams by linearity.

We will use $Z$ and $\fw_\fg$ to compute $J$ in Equation~(\ref{rep}).
See also \rmkref{rmk_color}.

\subsubsection{A lemma on $J$}\label{lemma}
The symmetric algebra $\fS(\Pi_\fh)$ in the simple roots $\Pi_\fh$ can be considered 
as polynomials on the root lattice as follows.
For any monomial $p = \ga_{i_1}\cdots\ga_{i_n}$ in $\fS(\Pi_\fh)$ and any 
$\mu\in Y$, let 
$$
p(\gm) = (\ga_{i_1} | \gm)\cdots(\ga_{i_n} | \gm)
$$
(cf.\ \cite{le3}). The degree of $p$ is $n$.

Recall from \secref{sec:tau} that 
$J = J^\fg$ is the quantum invariant of $\cV_2$-colored ribbon graphs in 
$S^3$ defined in I.2.5 of \cite{tu1}. Suppose that $\Omega$ is a $\cV_2$-colored ribbon graph is $S^3$.
Let $L$ be an $m$-component link in $S^3$ disjoint from $\Omega$. 
Let $\mu = (\mu_1,\ldots,\mu_m)$ with $\mu_i\in Y\cap C_k$. 
Denote by $L(\gm)$ the $\cV_2$-colored framed link with the $i$-th component colored by $\tVmi$.
Therefore $\Omega\sqcup L(\mu)$ is a $\cV_2$-extended ribbon graph.
The following lemma deals with the dependence
of $J_{\gW\sqcup L(\mu)}$ 
on $\gm_1, \ldots, \gm_m$.
Again, we assume that $r$ is an odd prime $\ge d h^\vee$ which 
does not divide $\det(a_{ij})$ and $\xi$ is a primitive $r$-th root of 1.

\begin{lemma}\label{lemma_new}
In the above notation,
\be\label{eqn_prop}
J_{\gW\sqcup L(\mu)} = \sum_{i=0}^{m(r\ell - \dim \fg)/2-1}
f_i(\gm_1, \ldots, \gm_m)\,(v-1)^i|_{v=\xi} + R
\ee
where $f_i$ is a polynomial in $\gm_1, \cdots, \gm_m$ 
of total degree not exceeding $2 i + m s$ for all $i$ and $R\in\zx$ is divisible by 
$(\xi -1)^{m(r\ell - \dim \fg)/2}$ in $\zx$.
Here $s$ is the number of
positive roots of $\fg$. Furthermore, for any $i$, the polynomial $f_i$
takes integer values for $\gm_1,\ldots, \gm_m\in C_k\cap Y$. 
\end{lemma}
\begin{proof}
It's obvious that $\gW\sqcup L$ 
can be presented as a diagram shown in \figref{fig}, where $T_2$ is the diagram of a string tangle and $T_1$ is the diagram of a ribbon graph such that all coupons of $\gW$ are inside $T_1$ 
and the link $L$ is inside $T_2$. To see this we first draw an arbitrary diagram presenting $\gW\sqcup L$ inside a disk
$D_1$. Then we pull $L$ into a disk $D_2$ disjoint from $D_1$ keeping all coupons inside $D_1$. After that 
the strings connecting $D_1$ and $D_2$, which are disjoint from $L$, can be arranged
into the diagram shown in \figref{fig}.
\begin{figure}\anchor{fig}\small
\centering
\psfrag{1}{$T_1$}
\psfrag{2}{$T_2$}
\psfrag{4}{$U_1$}
\psfrag{3}{$U_n$}
\includegraphics[width=3in]{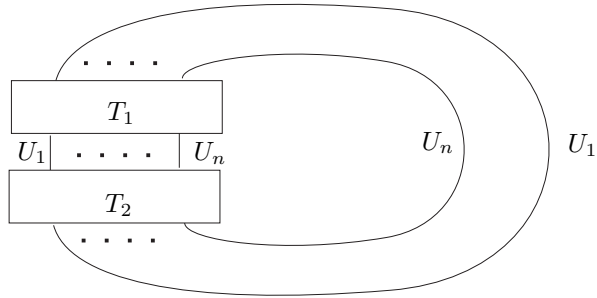}
\caption{A tangle presentation for $\gW\sqcup L$}\label{fig}
\end{figure}
Hence 
$$
J_{\gW\sqcup L(\mu)} =
\tr_q\left(J_{T_1}\circ J_{\Ttm} \right) = \tr\left(J_{T_1}\circ J_{\Ttm}\circ K_{2\rho}\right).
$$
Suppose that the strings connecting $T_1$ and $T_2$ have colors $U_1,\ldots, U_n$. Then
both $J_{T_1}$ and $J_{\Ttm}$ are endomorphisms of 
$U = \otimes_iU_i^{\epsilon_i}$ where $\epsilon_i=\pm$ is determined
by the orientation of the $i$-th string, cf.\ \secref{rg}.

Recall from \secref{4cat} that for any object $V$ in $\cV_2$ there exists a unique object $V'$ in
$\cV_1$ such that $V=\bar{V'}$. Pick a
basis $B'$ for $U'$ which induces a basis $B$ for $U$.
Then $J_{T_1}$ is a matrix over $\zx$ with respect
to $B$. Let $\hJ_{T_1}$ be the matrix over $\cA$ obtained from $J_{T_1}$ by substituting $\xi$ by $v$.
Note that $\hJ_{T_1}$
is not well-defined. For example $\xi^r$ can be substituted by either $v^r$ or 1. But one certainly has
$\hJ_{T_1}|_{v=\xi}=J_{T_1}$.

Let $\Ttmp$ be a $\cV_1$-colored framed string tangle obtained
from $T_2(\mu)$ by changing color on each of its component, including the bands, from $V$ to $V'$. 
Note that $T_2$ contains no coupons.
Let $J'$ be the invariant for $\cV_1$-colored
ribbon graphs constructed in Theorem I.2.5 of \cite{tu1}.
Then $J'_{\Ttmp}: U' \to U'$ is a morphism in $\cV_1$ which
is presented as a matrix
over $\cA$ with respect to $B'$. Since the ribbon structures on $\cV_1$ and $\cV_2$
are induced by the same quasi-$R$-matrix $\Theta$ (cf.\ Equations~(\ref{eq:quasi-R}-\ref{eq:Twist})), $J'_{\Ttmp}|_{v=\xi}$ is the matrix over $\zx$
presenting $J_{T_2(\mu)}$ with respect to $B$.
We have
\be\label{new}
J_{\gW\sqcup L(\mu)} = \tr\left(\hJ_{T_1} J'_{\Ttmp} K_{2\rho}\right)|_{v=\xi},
\ee
where $\hJ_{T_1} J'_{\Ttmp} K_{2\rho}$ is considered as a product of matrices.
By \lemmaref{last} below we have 
\be\label{in_h}
\tr\left(\hJ_{T_1} J'_{\Ttmp} K_{2\rho}\right)|_{v=e^{h/2}} = \sum_{i=0}^{\infty}
p_i(\gm_1, \ldots, \gm_m)\,h^i
\ee
where $p_i$ is a polynomial in root lattice of degree at most
$2 i + m s$ for all $i$. (Note that $K_{2\rho}$ and $\hJ_{T_1}$ are independent of $\mu$.) Substituting $v$ back into Equation~(\ref{in_h}), i.e.\ letting
$$
h = 2\ln(v) = 2\sum_{i=1}^\infty (-1)^{i-1} (v-1)^i/i,
$$
we have

\begin{eqnarray*}
\tr\left(\hJ_{T_1} J'_{\Ttmp} K_{2\rho}\right) &=&\sum_{i=0}^\infty f_i(\gm_1,\ldots,\gm_m) (v-1)^i \\
&=& \sum_{i=0}^{m(r\ell - \dim \fg)/2-1}
f_i(\gm_1, \ldots, \gm_m)\,(v-1)^i + R_v
\end{eqnarray*}
with $R_v$ divisible by $(v -1)^{m(r\ell - \dim \fg)/2}$ in $\cA$. The element $R_v$ is {\it a priori} in $\BC[[v]]$. It belongs to $\cA$ follows from
the fact that $\hJ_{T_1} J'_{\Ttmp} K_{2\rho}$
can be presented as a matrix over $\cA$ with respect to basis $B'$ if $\mu_i\in Y_+$. For the same reason we see that $f_i$ takes integer values if $\mu_i \in Y_+$. 
The coefficients $f_i$ are polynomials in root lattice of degree at most $2 i + m s$.
Substituting $v$ by $\xi$
we get Equation~(\ref{eqn_prop}). 
\end{proof}

Before we state \lemmaref{last} let's recall a fact first.

\subsubsection{A fact about $\Ug$}\label{Ug}
Any element $x$ in $\Ug$ can be written uniquely as a sum
$
x = z + u
$ 
for some $z$ in $\cen(\Ug)$ and $u$ in $[\Ug, \Ug]$
such that $\deg(z)\leq \deg(x)$ (cf.\ Page 105 Exercise 19 of \cite{bourbaki} or
Section 2.3.3 of \cite{warner}\footnote{We thank the referee for pointing out this reference to us.}).
Here $\deg$ is the PBW-degree.
Let $V_\lambda$ be the simple $\fg$-module with highest weight $\lambda\in Y_+$.
Then $z$ acts on $V_\lambda$ as a scalar $z(\lambda)\in\BC$. By the 
Harish-Chandra theory 
$z(\lambda)$ is a polynomial function on root lattice with degree no greater
than $\deg(z)$. Then by the Weyl formula the trace 
$\tr_{V_\lambda}(z) = \dim(V_\lambda) z(\lambda)$ is a polynomial on root
lattice with degree at most $\deg(z) + s$. Since $u$ has trace 0 it follows that $\tr_{V_\lambda}(x)$
is also a polynomial in root lattice which has degree at most $\deg(x) + s$.

\subsubsection{A lemma on $J'$}\label{J'}
Recall that $J'$ is the invariant for $\cV_1$-colored
ribbon graphs constructed in Theorem I.2.5 of \cite{tu1}. 
The following lemma is a generalization of Lemma 4.6 of \cite{le2}.

\begin{lemma}\label{last}
For $\nu=(\nu_1,\ldots,\nu_m)$, let $T(\nu)$ be a $\cV_1$-colored
framed string tangle whose closed components are colored by ${_1\! V_{\nu_1}},\ldots,{_1\! V_{\nu_m}}$.
Then, $J'_{T(\nu)}$ can be presented by a matrix $(c_{ij})$ such that 
\be\label{cij}
c_{ij}|_{v=e^{h/2}} = \sum_{k=0}^\infty p_k(\mu_1,\ldots,\mu_m)h^k
\ee
where $p_k$ is a polynomial in $\mu_1,\ldots,\mu_m$ of total degree not exceeding $2i+ms$.
\end{lemma}

Let $\Vg$ be the category of finite
dimensional $\fg$-modules.
It is well known that for each object $V$ in $\cV_1$ there exists
a unique object $V_\fg$ in $\Vg$ with the same
weight space as $V$.

\begin{proof}
Let $T(\tilde\nu)$ be the $\Vg$-colored framed string tangle 
obtained from $T(\nu)$ by changing all colors from $V$ to $V_\fg$.
Then one can calculate $J'_{T(\nu)}$ from the Kontsevich integral (Equation~(\ref{eq:kon}))
$$
Z(T(\tilde\nu)) = \sum_{i=0}^\infty Z_i(T(\tilde\nu)) h^i
$$
as follows. Suppose that $J'_{T(\nu)}$ is an endomorphism of
$V\in \cV_1$. Then it follows that $\fw_\fg(Z_i(T(\tilde\nu)))$ is an endomorphism of
$V_\fg\in \Vg$, cf.\ Equation~(\ref{eq:we}).
For any basis of $V$ and of $V_\fg$ there exists an invertible matrix $A$ over $\BC[[h]]$, independent of $\nu$, 
such that (as matrices over $\BC[[h]]$) we have
\be\label{rep}
J'_{T(\nu)}|_{v=e^{h/2}}= A\left(\sum_{i=0}^\infty (-1)^i \fw_\fg(Z_i(T(\tilde\nu))) h^i\right)A^{-1}.
\ee
(See \rmkref{rmk_color} below.)

Hence each entry of
$J'_{T(\nu)}|_{v=e^{h/2}}$ is of the form
$$
\sum_{k=0}^\infty p_k(\nu_1,\ldots,\nu_m)h^k
$$
where $p_k$ is a polynomial in $\nu_1,\ldots,\nu_m$ of total degree not exceeding $2i+ms$ (each chord has two chord-ends).
Now Equation~(\ref{cij}) follows from this and (\ref{rep}).
\end{proof}

\begin{rmk}\label{rmk_color}
Equation~(\ref{rep}) is proved in Theorem 7.2 of \cite{kt}. Strictly speaking, one needs to generalize the Kontsevich integral $Z$ in \secref{kon} (resp.\ weight system $\fw_\fg$ in \secref{weight}) to the case when components of tangles (resp.\ chord diagrams) are colored by (not necessarily identical) objects in $\Vg$. 
But this generalization is quite straightforward.
\end{rmk}

\subsubsection{Proof of \propref{prop_poly}}\label{proof_poly}
Since $J$ is defined over $\cV_2$, $F_{(L, \gW)}\in\zx$.
By \lemmaref{lemma_new}
$$
J_{U^m(\mu)} J_{\gW\sqcup L(\mu)} = \sum_{i=0}^{m(r\ell - \dim \fg)/2-1}
f'_i(\gm_1, \ldots, \gm_m)\,(\xi-1)^i + R'
$$
where $R'\in\zx$ is divisible by $(\xi-1)^{m(r\ell - \dim \fg)/2}$ 
and $f'_i$'s are polynomials in root lattice of degree at most $2 i + 2 m s$ taking integer
values for $\mu_i \in Y_+\cap C_k$.
Since
$$
F_{(L, \gW)} = \sum_{\gm\in (Y_+\cap C_k)^m} J_{U^m(\mu)} J_{\gW\sqcup L(\mu)}
$$
we can apply Lemma 4.7 in \cite{le2}\footnote{The proof of this lemma refers to Corollary 4.14 of \cite{le3}, which
should be Corollary 4.11 instead.} to show that
$F_{(L, \gW)}$ is divisible by $(\xi - 1)^{m(r\ell - \dim \fg)/2}$ in $\zx$.
This ends the proof of \propref{prop_poly}.

\Addresses\recd

\end{document}

%% file: agtout.tex
\def\ifplaintex{\expandafter\ifx\csname documentclass\endcsname\relax}

\def\gtp{{\mathsurround=0pt\it $\cal G\mskip-2mu$eometry \&\ 
$\cal T\!\!$opology $\cal P\!$ublications}}  

\def\recd{{\small Received:\qua\receiveddate\ifx\reviseddate\relax
\else\qquad Revised:\qua\reviseddate\fi\par}} 

\def\lognumber#1{\def\thelognumber{#1}}
\def\volumenumber#1{\def\thevolumenumber{#1}}
\def\volumeyear#1{\def\thevolumeyear{#1}}
\def\papernumber#1{\def\thepapernumber{#1}}
\def\pagenumbers#1#2{\def\startpage{#1}\def\finishpage{#2}}
\def\published#1{\def\publishdate{#1}}

\def\received#1{\def\receiveddate{#1}}
\def\revised#1{\def\reviseddate{#1}}
\def\accepted#1{\def\accepteddate{#1}}

\def\asciiaddress#1{\def\theasciiaddress{#1}}
\def\asciiemail#1{\def\theasciiemail{#1}}

\long\def\asciiabstract#1{\long\def\theasciiabstract{#1}}

\let\\\par\let\thelognumber\relax\let\thevolumenumber\relax
\let\thepapernumber\relax\let\thevolumeyear\relax\let\startpage\relax
\let\finishpage\relax\let\publishdate\relax\let\receiveddate\relax
\let\reviseddate\relax\let\accepteddate\relax\let\theasciititle\relax
\let\theasciiauthors\relax\let\theasciiaddress\relax
\let\theasciiabstract\relax

\let\theasciiemail\relax

\ifplaintex
\font\logobig=cmssbx10 scaled 3836
\font\logomed=cmssbx10 scaled 2557
\else
\font\logobig=cmssbx10 scaled 4200
\font\logomed=cmssbx10 scaled 2800
\fi

\long\def\makeagttitle{   
\count0=\startpage
\agt\hfill      
\hbox to 45truept{\vbox to 0pt{\vglue -13truept{\logomed A\kern -.37em{\logobig 
T}\kern -.38em G}\vss}\hss}
\break
{\small Volume \thevolumenumber\ (\thevolumeyear)
\startpage--\finishpage\nl
Published: \publishdate}

\vglue .25truein

{\parskip=0pt\leftskip 0pt plus
1fil\def\\{\par\smallskip}{\Large\bf\thetitle}\par\medskip} \vglue
0.05truein

{\parskip=0pt\leftskip 0pt plus 1fil\def\\{\par}{\sc\theauthors}
\par\medskip}%
 
\vglue 0.03truein 

{\small\leftskip 25truept\rightskip 25truept{\bf Abstract}\stdspace\theabstract

{\bf AMS Classification}\stdspace\theprimaryclass
\ifx\thesecondaryclass\relax\else; \thesecondaryclass\fi\par
{\bf Keywords}\stdspace \thekeywords\par}\vglue 7truept

}   

\ifplaintex
\font\phead=cmsl9 scaled 950
\font\pnum=cmbx10 scaled 913
\font\pfoot=cmsl9 scaled 950
\headline{\vbox to 0pt{\vskip -4.5mm\line{\small\phead\ifnum
\count0=\startpage ISSN 1472-2739 (on-line) 1472-2747 (printed)
\hfill {\pnum\folio}\else\ifodd\count0\def\\{ }%
\ifx\theshorttitle\relax\thetitle\else\theshorttitle\fi\hfill{\pnum\folio}
\else\def\\{ and }{\pnum\folio}\hfill\ifx\theshortauthors\relax\theauthors
\else\theshortauthors\fi\fi\fi}\vss}}
\footline{\vbox to 0pt{\vglue 0mm\line{\small\pfoot\ifnum\count0=\startpage
\copyright\ \gtp\hfill\else
\agt, Volume \thevolumenumber\ (\thevolumeyear)\hfill\fi}\vss}}
\else
\font=cmsl9 scaled 1050
\font\lnum=cmbx10 
\font=cmsl9 scaled 1050
\makeatletter
\def\@oddhead{{\small\lhead\ifnum\count0=\startpage ISSN 1472-2739 
(on-line) 1472-2747 (printed)\hfill {\lnum\number\count0}\else\ifodd\count0
\def\\{ }\ifx\theshorttitle\relax \thetitle \else\theshorttitle\fi\hfill
{\lnum\number\count0}\else\def\\{ and }{\lnum\number\count0}
\hfill\ifx\theshortauthors\relax 
\theauthors\else\theshortauthors\fi\fi\fi}}\def\@evenhead{\@oddhead}
\def\@oddfoot{\small\lfoot\ifnum\count0=\startpage\copyright\ \gtp\hfill\else
\agt, Volume \thevolumenumber\ (\thevolumeyear)\hfill\fi}
\def\@evenfoot{\@oddfoot}
\makeatother
\fi
\let\maketitlepage\makeagttitle
\let\makeshorttitle\maketitlepage
\let\maketitle\maketitlepage

\newwrite\gtoutfile
\long\gdef\makeheadfile{  
{\def\\{, }\def\s{ }
\immediate\openout\gtoutfile head.xxx
\immediate\write\gtoutfile{Proxy-for: \ifx\theasciiauthors\relax
\theauthors\else\theasciiauthors\fi\s<\ifx\theasciiemail\relax\theemail\else\theasciiemail\fi>}
\immediate\write\gtoutfile{\noexpand\\}
\immediate\write\gtoutfile{Authors: \ifx\theasciiauthors\relax
\theauthors\else\theasciiauthors\fi}
{\def\\{ }\immediate\write\gtoutfile{Title: \ifx\theasciititle\relax
\thetitle\else\theasciititle\fi}}
\immediate\write\gtoutfile{Subj-class: GT or SG, GR etc}
\immediate\write\gtoutfile{MSC-class: \theprimaryclass\ifx\thesecondaryclass\relax\else, \thesecondaryclass\fi}
\immediate\write\gtoutfile{Journal-ref: Algebr. Geom. Topol. \thevolumenumber\s
(\thevolumeyear) \startpage-\finishpage}
\immediate\write\gtoutfile{Comments: Published by Algebraic and
Geometric Topology at}
\immediate\write\gtoutfile{\s\s\s  http://www.maths.warwick.ac.uk/agt/AGTVol\thevolumenumber/agt-\thevolumenumber-\thepapernumber.abs.html}
\immediate\write\gtoutfile{\noexpand\\}
\immediate\write\gtoutfile{}
\ifx\theasciiabstract\relax
\immediate\write\gtoutfile{\theabstract}\else
\immediate\write\gtoutfile{\theasciiabstract}\fi
\immediate\write\gtoutfile{}
\immediate\write\gtoutfile{\noexpand\\}
\immediate\write\gtoutfile{}
\immediate\closeout\gtoutfile}}  

\def\maketitlepage{\makeagttitle\makeheadfile}
\let\makeshorttitle\maketitlepage
\let\maketitle\maketitlepage

%% file: agt-5-50.bbl
\begin{thebibliography}

\bibitem[BK]{bk}
\textbf{B Bakalov}, \textbf{A Kirillov, Jr}, \emph{Lectures on tensor
  categories and modular functors}, University Lecture Series 21, American
  Mathematical Society, Providence, RI (2001) \MR{1797619}

\bibitem[BHMV]{bhmv}
\textbf{C Blanchet}, \textbf{N Habegger}, \textbf{G Masbaum}, \textbf{P Vogel},
  \emph{Topological quantum field theories derived from the {K}auffman
  bracket}, Topology 34 (1995) 883--927 \MR{1362791}

\bibitem[B]{bourbaki}
\textbf{N Bourbaki}, \emph{Lie groups and {L}ie algebras. {C}hapters 1--3},
  Elements of Mathematics (Berlin), Springer-Verlag, Berlin (1989) \MR{979493}

\bibitem[C]{chen}
{\bf Qi Chen},
 {\it On certain integral tensor categories and integral TQFTs},
 arXiv:math.QA/0408356

\bibitem[CL]{chenle1}
\textbf{Q Chen}, \textbf{T\,T\,Q Le}, \emph{Quantum invariants of periodic links
  and periodic 3-manifolds}, Fund. Math. 184 (2004) 55--71 \MR{2128042}

\bibitem[GKP]{gkp}
\textbf{P\,M Gilmer}, \textbf{J Kania-Bartoszynska}, \textbf{J\,H Przytycki},
  \emph{3-manifold invariants and periodicity of homology spheres}, 
\href{http://www.maths.warwick.ac.uk/agt/AGTVol2/agt-2-34.abs.html}%
{Alg. Geom. Topol. 2 (2002) 825--842} \MR{1936972}

\bibitem[G]{g2}
\textbf{P\,M Gilmer}, \emph{Integrality for {TQFT}s}, Duke Math. J. 125 (2004)
  389--413 \MR{MR2096678}

\bibitem[K]{kassel}
\textbf{C Kassel}, \emph{Quantum groups}, Graduate Texts in Mathematics 155,
  Springer-Verlag, New York (1995) \MR{1321145}

\bibitem[KT]{kt}
\textbf{C Kassel}, \textbf{V Turaev}, \emph{Chord diagram invariants of tangles
  and graphs}, Duke Math. J. 92 (1998) 497--552 \MR{1620522}

\bibitem[K]{kirillov}
\textbf{A\,A Kirillov, Jr}, \emph{On an inner product in modular tensor
  categories}, J. Amer. Math. Soc. 9 (1996) 1135--1169 \MR{1358983}

\bibitem[L1]{le3}
\textbf{T\,T\,Q Le}, \emph{On perturbative {${\rm PSU}(n)$} invariants of
  rational homology 3-spheres}, Topology 39 (2000) 813--849 \MR{1760430}

\bibitem[L2]{le2}
\textbf{T\,T\,Q Le}, \emph{Quantum invariants of 3-manifolds: integrality,
  splitting, and perturbative expansion}, Topology Appl. 127 (2003) 125--152
  \MR{1953323}

\bibitem[LM]{lm}
\textbf{T\,Q\,T Le}, \textbf{J Murakami}, \emph{The universal
  {V}assiliev-{K}ontsevich invariant for framed oriented links}, Compositio
  Math. 102 (1996) 41--64 \MR{1394520}

\bibitem[L]{lu}
\textbf{G Lusztig}, \emph{Introduction to quantum groups}, Progress in
  Mathematics 110, Birkh\"auser Boston Inc. Boston, MA (1993) \MR{1227098}

\bibitem[RT]{rt}
\textbf{N Reshetikhin}, \textbf{V\,G Turaev}, \emph{Invariants of
  {$3$}-manifolds via link polynomials and quantum groups}, Invent. Math. 103
  (1991) 547--597 \MR{1091619}

\bibitem[T]{tu1}
\textbf{V\,G Turaev}, \emph{Quantum invariants of knots and 3-manifolds}, de
  Gruyter Studies in Mathematics 18, Walter de Gruyter \& Co. Berlin (1994)
  \MR{1292673}

\bibitem[W]{warner}
{\bf Garth Warner},
 {\em Harmonic analysis on semi-simple {L}ie groups, {I}},
 Grundlehren series 188, Springer-Verlag (1972) \MR{0498999}


\end{thebibliography}
